\documentclass[a4paper,11pt]{amsart}
\usepackage{amssymb}
\usepackage{amsthm}
\usepackage{mathrsfs}
\usepackage{latexsym}
\theoremstyle{plain}
\newtheorem{Theorem}{Theorem}
\newtheorem{Corollary}{Corollary}

\newtheorem{Conjecture}{Conjecture}
\theoremstyle{definition}

\newtheorem{Remark}{Remark}

\newcommand{\vdim}{\mathrm{vdim}}
\newcommand{\Pic}{\mathrm{Pic}}
\newcommand{\Hilb}{\mathrm{Hilb}}

\newcommand{\R}{\mathbb{R}}
\newcommand{\Q}{\mathbb{Q}}
\newcommand{\Z}{\mathbb{Z}}

\newcommand{\mgnbar}{\overline{\mathcal{M}}_{g,n}}
\newcommand{\monbar}{\overline{\mathcal{M}}_{0,n}}
\newcommand{\BB}{\mathscr B}

\title[Quantum cohomology of moduli of genus zero curves]{Quantum cohomology of moduli spaces \\ of genus zero stable curves}

\author{Claudio Fontanari}

\email{claudio.fontanari@polito.it}\curraddr{
Dipartimento di Matematica \\ Politecnico di Torino \\
Corso Duca degli Abruzzi 24 \\ 10129 Torino \\ Italy.}

\date{}


\begin{document}

\begin{abstract}
We investigate the (small) quantum cohomology ring of the moduli spaces 
$\monbar$ of stable $n$-pointed curves of genus $0$. In particular, we 
determine an explicit presentation in the case $n=5$ and we outline a 
computational approach to the case $n=6$.
\end{abstract}

\maketitle

\section{Introduction} 

The (small) quantum cohomology ring of a smooth algebraic variety with 
$\star$-product defined in terms of ($3$-point) Gromov-Witten invariants 
is a formal deformation of the classical Chow ring in the sense that 
the $\star$-product specializes to the cup-product when the formal 
parameters are set to $0$. 
The notion of quantum Chow ring has been recently extended also to 
smooth orbifolds and its degree zero part is usually called the 
stringy Chow ring. 
In \cite{AGV:02} the stringy Chow ring of $\overline{\mathcal{M}}_{1,1}$ 
has been computed, while in \cite{Spencer:06} the case of $\mathcal{M}_2$ 
is handled and that of $\overline{\mathcal{M}}_{2}$ is announced. 

Here instead we address the (small) quantum cohomology of the moduli 
spaces $\monbar$ of stable $n$-pointed curves of genus $0$. Even though 
these spaces are smooth projective varieties with a quite explicit 
description, nonetheless their geometry turns out to be rather involved: 
indeed, just to quote a couple of astonishing facts, we mention that, 
despite the serious efforts by many valuable mathematicians, their 
ample (resp., effective) cone has been determined so far only for $n \le 7$
and $n \le 6$ resp. (see \cite{KMK:96} and \cite{HT:02} resp.). 

In the present paper we provide an explicit presentation of the small 
quantum cohomology ring in the case $n=5$ (see Corollary~\ref{n=5}) 
relying on previous work by G\"ottsche and Pandharipande (\cite{GP:98}) 
and we suggest a computational approach to the case $n=6$ (see 
Remark~\ref{n=6}) inspired by Gathmann \cite{G:01} (see also \cite{BM:04}, where the small quantum cohomology of all del Pezzo surfaces is calculated, and \cite{B:04}, where a general theorem about semisimplicity conservation under blowing-up of points is proved). 

The author is grateful to Gianfranco Casnati, Gianni Ciolli, and Barbara Fantechi for inspiring conversations and enlightening suggestions. Thanks are also due to 
Yuri I. Manin and Dan Abramovich for kindly pointing 
out references \cite{B:04}, \cite{BM:04}, and 
\cite{BK:05}, respectively.

\section{Preliminaries} 

\subsection{(Small) quantum cohomology}\label{quantum}

Let $X$ be a smooth complex projective variety, let $T_0 = 1 \in A^0(X), 
T_1, \ldots,$ $T_m$ be a homogeneous basis of the graded vector space 
$V := H^*(X, \Q)$. Let $T^0 = \mathrm{point}, T^1, \ldots, T^m \in V$ 
be the (Poincar\'e) dual basis and let $E \subset H_2(X, \Z)$ denote 
the subset of effective curves. 

The (small) quantum cohomology ring of $X$ is a $\star$-product structure 
on $V \otimes R$, where $R$ is a formal power series ring and the quantum 
product $\star$ reduces to the usual cup product $\cup$ when all formal 
variables are set to zero. Namely, given classes $\alpha_1$ and $\alpha_2
\in H^*(X, \Q)$, their quantum product is defined as follows:
$$
\alpha_1 \star \alpha_2 = \alpha_1 \cup \alpha_2 + \sum_{\beta \in E} 
\sum_{i=0}^m I_\beta(\alpha_1, \alpha_2, T_i) T^i q^{\beta}
$$
where $I_\beta(\cdot,\cdot,\cdot)$ denotes the ($3$-point) Gromov-Witten 
invariant relative to class $\beta$ (see for instance \cite{FP:97}, \S~7). 

\subsection{Moduli spaces of genus zero (stable) curves}\label{moduli} 

The moduli spaces $\monbar$ of stable $n$-pointed curves of genus $0$ 
is a smooth projective variety of dimension $n-3$ which can be explicitely 
obtained from $\mathbb{P}^{n-3}$ via the following construction due to 
Kapranov. For every $n \ge 4$, let $X \subset \mathbb{P}^{n-3}$ be a set of $n-1$ 
points in linear general position, let $B^0 := \mathbb{P}^{n-3}$ and for $i \ge 0$ 
let $B^{i+1} \to B^i$ be the blow-up of $B^i$ along the proper transforms of 
the $\binom{n-1}{i+1}$ $i$-planes through $i+1$ points. With the above 
notation, we have $\monbar \cong B^{n-4}$ (see for instance \cite{V:02}).

Let $P:= \{1,2, \ldots, n \}$ and for every $S \subset P$ with 
$2 \le \vert S \vert \le n-2$ let $\Delta_{\{0,S\}}$ be the boundary 
component of $\monbar$ whose general element is the union of 
two copies of $\mathbb{P}^1$, labelled respectively by $S$ 
and $P \setminus S$, meeting at one point. We denote by 
$\delta_S$ the corresponding class in $\Pic(\monbar)$ and    
we define inductively:
\begin{eqnarray*}
\BB_4 &:=& \{ \delta_{\{2,3\}} \}\\
\BB_i &:=& \BB_{i-1} \cup \{ \delta_B: B \subseteq \{1, \ldots, i \}, 
i \notin B \supseteq \{i-1, i-2 \} \}\\ 
& &\cup \{\delta_{B^c \setminus \{ i \}}: \delta_B \in 
\BB_{i-1} \setminus \BB_{i-2} \}.
\end{eqnarray*}
Then according to \cite{F:05}, Proposition~1, $\BB_n$ is a basis of $\Pic(\monbar)$. 

\section{The results}

\subsection{The case $n=5$} 
In order to determine the quantum cohomology ring we first need to manage
the quantum product between two divisor classes.

\begin{Theorem}~\label{main}
Let $X$ be $\mathbb{P}^2$ blown up at $4$ points in linear general position. 
Let $H, E_1, \ldots, E_4$ denote the strict transform of the hyperplane class 
and the exceptional divisor classes respectively. 
If $\Delta_1, \Delta_2 \in H^2(X, \Q)$ then their quantum product can be expressed 
as follows:
\begin{eqnarray*}
\Delta_1 \star \Delta_2 &=& \Delta_1 . \Delta_2 
- \sum_{i=1}^4 \Delta_1.E_i 
\Delta_2.E_i E_i q^{0, e_i} \\
& & + \sum_{1 \le i < j \le 4} \Delta_1.(H-E_i-E_j)
\Delta_2.(H-E_i-E_j)\\
& &+ (H+E_i+E_j)q^{1, e_i+e_j}\\
& &+ \sum \Delta_1.(H-\sum_{i=1}^4 \varepsilon_i E_i) 
\Delta_2.(H-\sum_{i=1}^4 \varepsilon_iE_i) q^{1, \sum_{i=1}^4 \varepsilon_ie_i}\\
& & + \sum \Delta_1. (2H-\sum_{i=1}^4 \varepsilon_iE_i)
\Delta_2. (2H-\sum_{i=1}^4 \varepsilon_iE_i) 
q^{2,\sum_{i=1}^4 \varepsilon_ie_i}
\end{eqnarray*}
where the sums run over $\varepsilon_i \in \{0,1 \}$ and if 
$\beta = a H - \sum_{i=1}^4 b_i E_i$ we denote $q^{\beta}$ by
$q^{a,(b_1,b_2,b_3,b_4)}$ writing $e_i$ for the $i$-th vector of 
the canonical basis of $\R^4$.
\end{Theorem}
 
\proof By \cite{BP:04}, Corollary~3.3, the effective cone of $X$ is 
generated by the following divisors: 
\begin{eqnarray*}
D_1 &=& H - E_1 - E_2 \\
D_2 &=& H - E_1 - E_3 \\
D_3 &=& H - E_1 - E_4 \\
D_4 &=& H - E_2 - E_3 \\
D_5 &=& H - E_2 - E_4 \\
D_6 &=& H - E_3 - E_4 \\
D_7 &=& E_1 \\
D_8 &=& E_2 \\
D_9 &=& E_3 \\
D_{10} &=& E_4 \\
\end{eqnarray*}
Hence if 
$$
\beta = d H - \sum_{i=1}^4 a_i E_i  
$$
is an effective curve, then we can write 
$$
\beta = \sum_{i=1}^{10} c_i D_i 
$$
with $c_i \ge 0$ for every $i$, in particular we have 
$$
d = \sum_{i=1}^6 c_i
$$
By the divisor axiom, 
$$
I_\beta(\Delta_1,\Delta_2,T_i) = (\Delta_1.\beta)(\Delta_2.\beta) I_\beta(T_i)
$$
therefore we need only to compute $1$-point Gromov-Witten invariants. 
Since 
$$
-K_X = \mathcal{O}_X(1) = 3H - E_1 - E_2 - E_3 - E_4
$$
we deduce that the expected dimension of the corresponding moduli space of stable 
maps is
\begin{equation}\label{dimension}
\vdim \overline{\mathcal{M}}_{0,1}(X, \beta) = - K_X. \beta = \sum_{i=1}^{10} c_i 
\end{equation}
We have to consider separately the three cases $T_i = 1$, $T_i = D$ a divisor 
and $T_i = \mathrm{point}$. It is a general fact that $I_\beta(1)=0$ (see for 
instance \cite{FP:97}, \S~7.(II)). Next, if $D$ is a divisor then 
$I_\beta(D) \ne 0$ only if $\vdim \overline{\mathcal{M}}_{0,1}(X, \beta)=1$, 
in particular we have $d \le 1$. If $d=0$, then $\beta$ is purely exceptional 
and from \cite{G:01}, Lemma~2.3~(i), it follows that $I_\beta(D) \ne 0$ if 
and only if $\beta = E_i$, $D = E_i$ and $I_{E_i}(E_i)=-1$. If instead $d=1$, 
then $\beta = H - E_i - E_j$ and $I_\beta(D) \ne 0$ if and only if $D$ is either 
$H$, or $E_i$, or $E_j$, and $I_\beta(D)=1$. Finally, if $T_i = \mathrm{point}$ 
then $I_\beta(\mathrm{point}) \ne 0$ only if 
$\vdim \overline{\mathcal{M}}_{0,1}(X, \beta)=2$, in particular we have $d \le 2$.
If $d=0$, then $\beta$ is purely exceptional and $I_\beta(\mathrm{point}) = 0$
by \cite{G:01}, Lemma~2.3~(i). If $d=1$, we have $I_\beta(\mathrm{point}) = 1$
for $\beta = H - \sum_{i=1}^4 \varepsilon_i E_i$ with $\varepsilon_i \in \{0,1 \}$ 
and $I_\beta(\mathrm{point}) = 0$ otherwise by \cite{GP:98}, 
\S~5.2 and \S~3.(P4)--(P5). If $d=2$, we have $I_\beta(\mathrm{point}) = 1$
for $\beta = 2H - \sum_{i=1}^4 \varepsilon_i E_i$  with $\varepsilon_i \in \{0,1 \}$ 
and $I_\beta(\mathrm{point}) = 0$ otherwise by \cite{GP:98}, \S~5.2
and \S~3.(P4)--(P5). Hence our claim follows. 

\qed

As a consequence, we can perform the computation we are interested in. 

\begin{Corollary}\label{n=5}
The small quantum cohomology ring of $\overline{\mathcal{M}}_{0,5}$ 
admits the following explicit presentation:
$$
QH^*_s(\overline{\mathcal{M}}_{0,5}) = \frac{\Q[q^{0, \sum_{i=1}^4 \varepsilon_i e_i}, 
q^{1,\sum_{i=1}^4 \varepsilon_i e_i}, q^{2, \sum_{i=1}^4 \varepsilon_i e_i},
\delta_{2,3}, \delta_{3,4},\delta_{1,5},
\delta_{2,5},\delta_{1,4}]}{(f_i^*)_{i=1, \ldots, 5}}  
$$
where $\varepsilon_i \in \{0,1\}$, $e_i$ denotes the $i$-th vector of 
the canonical basis of $\R^4$ and
\begin{eqnarray*}
f_1^* &=& \delta_{2,3} \star \delta_{3,4}+ E_1 q^{0,(1,0,0,0)}-q^{1,(0,0,0,0)}
-q^{1,(0,0,1,0)}-q^{1,(1,1,0,1)}\\
& &-q^{1,(1,1,1,1)}-4q^{2,(0,0,0,0)}-q^{2,(1,0,0,0)}-2q^{2,(0,1,0,0)}-4q^{2,(0,0,1,0)}\\
& &-2q^{2,(0,0,0,1)}-q^{2,(1,0,1,0)}-2q^{2,(0,1,1,0)}-q^{2,(0,1,0,1)}\\
& &-2q^{2,(0,0,1,1)}-q^{2,(0,1,1,1)}\\
f_2^* &=& \delta_{2,3} \star \delta_{2,5}-(H+E_2+E_3)q^{1,(0,1,1,0)}-q^{1,(0,1,0,0)}
-q^{1,(0,1,1,0)}\\
& &+q^{1,(1,1,0,1)}+q^{1,(1,1,1,1)}-2q^{2,(0,1,0,0)}-q^{2,(1,1,0,0)}
-2q^{2,(0,1,1,0)}\\
& &-q^{2,(0,1,0,1)}-q^{2,(0,1,1,1)}
-q^{2,(1,1,1,0)}\\
f_3^* &=& \delta_{3,4} \star \delta_{1,4}+ E_2 q^{0,(0,1,0,0)}-q^{1,(0,0,0,0)}
-q^{1,(0,0,0,1)}-q^{1,(1,1,1,0)}\\
& &-q^{1,(1,1,1,1)}-4q^{2,(0,0,0,0)}-2q^{2,(1,0,0,0)}-q^{2,(0,1,0,0)}\\
& &-2q^{2,(0,0,1,0)}-4q^{2,(0,0,0,1)}-q^{2,(1,0,1,0)}-2q^{2,(1,0,0,1)}-q^{2,(0,1,0,1)}\\
& &-2q^{2,(0,0,1,1)}-q^{2,(1,0,1,1)}\\
f_4^* &=& \delta_{1,5} \star \delta_{2,5}-(H+E_1+E_2)q^{1,(1,1,0,0)}
-q^{1,(1,1,0,0)}-q^{1,(1,1,1,0)}\\
& &-q^{1,(1,1,0,1)}-q^{1,(1,1,1,1)}-q^{2,(1,1,0,0)}-q^{2,(1,1,1,0)}
-q^{2,(1,1,0,1)}\\
& &-q^{2,(1,1,1,1)}\\ 
f_5^* &=& \delta_{1,5} \star \delta_{1,4}-(H+E_1+E_4)q^{1,(1,0,0,1)}-q^{1,(1,0,0,0)}
-q^{1,(1,0,0,1)}\\
& &+q^{1,(1,1,1,0)}+q^{1,(1,1,1,1)}-2q^{2,(1,0,0,0)}-2q^{2,(1,0,0,1)}
-q^{2,(1,1,0,0)}\\
& &-q^{2,(1,0,1,0)}-q^{2,(1,1,0,1)}-q^{2,(1,0,1,1)}
\end{eqnarray*}
\end{Corollary}

\proof From Keel's results in \cite{K:92} we deduce the following presentation 
of the classical Chow ring of $\overline{\mathcal{M}}_{0,5}$ in terms of the 
basis $\BB_5$ recalled in \S~\ref{moduli}:
$$
A^*(\overline{\mathcal{M}}_{0,5}) = \frac{\Z[\delta_{2,3}, \delta_{3,4},\delta_{1,5},
\delta_{2,5},\delta_{1,4}]}{(f_i)_{i=1, \ldots, 5}}  
$$
where 
\begin{eqnarray*}
f_1^* &=& \delta_{2,3}.\delta_{3,4} = 0 \\
f_2^* &=& \delta_{2,3}.\delta_{2,5} = 0 \\ 
f_3^* &=& \delta_{3,4}.\delta_{1,4} = 0 \\
f_4^* &=& \delta_{1,5}.\delta_{2,5} = 0 \\ 
f_5^* &=& \delta_{1,5}.\delta_{1,4} = 0 
\end{eqnarray*}
According to Kapranov construction recalled in \S~\ref{moduli}, we can regard 
$\overline{\mathcal{M}}_{0,5}$ as $\mathbb{P}^2$ blown up at $4$ points in linear 
general position and obtain exactly as in \cite{V:02} the following identifications 
(here we take $5$ to be the special point):
\begin{eqnarray*}
\delta_{2,3} &=& H - E_1 - E_4 \\
\delta_{3,4} &=& H - E_1 - E_2 \\
\delta_{1,5} &=& E_1 \\
\delta_{2,5} &=& E_2 \\
\delta_{1,4} &=& H - E_2 - E_3
\end{eqnarray*}
Notice moreover that by (\ref{dimension}) $I_\beta(\alpha_1,\alpha_2,T_i) \ne 0$ 
only when $\sum c_i (-K_X.D_i)$ is a fixed number
with both $c_i \ge 0$ and $-K_X.D_i > 0$ for every $i$, hence there are only 
finitely many possible values for the exponents of the formal variables $q$ and 
the quantum cohomology ring turns out to be a polynomial ring.  
Hence our claim can be deduced from Theorem~\ref{main} by applying\cite{FP:97},
\S~10, Proposition~11 (see \cite{P:03}, Chapter~3, for analogous computations). 

\qed

\subsection{The case $n=6$}
Here we have obtained only partial results.

\begin{Conjecture}\label{associativity}
Let $Y$ be the blow-up a smooth projective threefold $X$ along a curve $C$ such that 
$g(C) \ge 1$ or $g(C)=0$ and $-K_X.C \ge 0$. Then the associativity equations of the 
quantum product suffice to determine all (genus $0$) Gromov-Witten invariants of $Y$ 
in terms of those of $X$.
\end{Conjecture}

\begin{Remark}
In order to address Conjecture~\ref{associativity}, one 
might wish to argue as in \cite{G:01}, proof of Theorem~2.1. 
Indeed, if both $\beta$ and $T$ are non-exceptional classes, 
then \cite{H:00}, Theorem~1.5, would even imply $I_\beta^Y(T) = I_\beta^X(T)$ (but see \cite{BK:05}, Remark~8, for a 
pertinent counterexample to the statement in \cite{H:00}). 
On the other hand, if $\beta$ is exceptional then 
$\beta = F$ and the only eventually nonzero invariants 
to be computed are $I_{dF}(d \varphi)$, where $F$ and $\varphi$ correspond to exceptional fibers. These invariants enumerate $d$-fold coverings of a fibre over a point in $C$, 
hence they are zero for $d \ge 2$ (otherwise a curve should lie in two different fibers). If instead $d=1$ then $I_F(\varphi)=-1$ (see for instance \cite{C:05}, Lemma~2). 
Unluckily, as far as we know, the analogue of \cite{G:01}, Algorithm~2.4, is still missing.  
\end{Remark}

\begin{Remark}\label{n=6}
From Conjecture~\ref{associativity} it would follow that all Gromov-Witten invariants 
of $\overline{\mathcal{M}}_{0,6}$ can be recursively computed. Indeed, as recalled in 
\S~\ref{moduli}, $\overline{\mathcal{M}}_{0,6}$ can be identified with $\mathbb{P}^3$ blown up in $5$ points in linear general position and along the cords between pairs of points. 
In order to check that $-K.C \ge 0$, let $C$ be the strict transform of the cord 
$l_{ab}$ between points $a$ and $b$ and choose planes $\pi$, $\rho$ in 
$\mathbb{P}^3$ such that $l_{ab} = \pi.\rho$. If $p: \tilde{X} \to X$  
denotes the blow up of $l_{ab}$ we have
\begin{eqnarray*}
p^* \pi &=& \tilde{\pi} + E_a + E_b \\
p^* \rho &=& \tilde{\rho} + E_a + E_b \\
K_{\tilde{X}} &=& K_{\mathbb{P}^3} + 2 \sum_i E_i + 2 \sum_{i,j} E_{ij} 
\end{eqnarray*}
(see \cite{H:77}, II., ex.~8.5) where $\tilde{\pi}$ and $\tilde{\rho}$ resp. are 
the strict transforms of $\pi$ and $\rho$ resp., while $E_i$ and $E_{ij}$ are 
the exceptional divisors corresponding to the points and the cords which have 
been previously blown up. Hence 
\begin{eqnarray*}
-K_{\tilde{X}}.C &=& (- K_{\mathbb{P}^3} - 2 \sum_i E_i - 2 \sum_{i,j} E_{ij}). \\ 
& & (p^* \pi - E_a - E_b) . (p^* \rho - E_a - E_b) \\
&=& - K_{\mathbb{P}^3}.\pi.\rho - 2 E_a^3 - 2 E_b^3 \\
&=& \mathcal{O}_{\mathbb{P}^3}(4).l_{ab}-2-2 = 0
\end{eqnarray*}
(recall that if $p: \tilde{X} \to X$ is the blow up of a smooth threefold along 
a smooth curve with exceptional divisor $E$ then $E.p^*C=0$ for every curve 
$C \subset X$, see for instance \cite{C:05}, Lemma~1). 
\end{Remark}

\end{document}